\documentclass[a4paper,10pt]{amsart}

\usepackage[utf8]{inputenc}
\usepackage[T1]{fontenc}
 \usepackage{amsfonts,mathrsfs}
 \usepackage{amsthm}
 \usepackage{amssymb}
 \usepackage{enumerate}
 \usepackage{tikz-cd}
 \usepackage{bm}
 \usepackage{csquotes}
 \usepackage{amsmath}
\usepackage{amssymb} 
 
\usepackage{pgfplots}
\usepackage{hyperref}
\hypersetup{colorlinks=true,linkcolor=black,citecolor=blue,urlcolor=blue}
\usepackage{booktabs}
\usepackage{array}
\usepackage{caption}
\usepackage{subcaption}

\theoremstyle{definition}

\theoremstyle{plain}

\usepackage[makeroom]{cancel}
\usepackage{graphicx}
\usepackage{cancel}
\usepackage{ulem}
\usepackage{setspace}

\usepackage{amsthm}

\usepackage{tikz}
\usepackage{verbatim}

\usepackage{amsrefs}

\usepackage{graphicx}
\usepackage{mathtools}

\usepackage{graphicx}
\usepackage{float} 
\usepackage{amsfonts}
\usepackage{amsthm}
\usepackage{latexsym}
\usepackage{amsmath}
\usepackage{amssymb}
\usepackage{amscd}
\usepackage{epsf}
\usepackage{tikz}
\usetikzlibrary{matrix,arrows}
\usepackage{enumerate}
\usepackage{eepic}
\usepackage{epic}
\usepackage{amsmath,amsbsy}
\usepackage{cancel}

\usepackage{amsthm}

\newtheorem*{conjecture*}{Conjecture}

\usepackage{amsrefs}

\theoremstyle{definition}

\theoremstyle{plain}
\newtheorem{theorem}{Theorem}
\newtheorem*{theorem*}{Theorem}

\usepackage{mathtools}

\title[Berry's phase under topology change]{Berry's phase under topology change}
\author{Pavel Kurasov}
\address{Stockholm University, Sweden} 
\email{kurasov@math.su.se}
\author{Vladislav Shubin}
\address{Stockholm University, Sweden} 
\email{v.v.choubine@gmail.com}
\author{Axel Tibbling}
\address{Stockholm University, Sweden}
\email{tibbling.axel@gmail.com}

\date{\today}

\begin{document}

\include{operators}
\include{fonts}


\nocite{*}

\begin{abstract}
Laplacians on metric graphs are used to construct continuous families of Hamiltonians with
different topological structure. One such family is used to demonstrate 
that Hamiltonians with
real-valued eigenfunctions may possess non-trivial geometric Berry's phase. 
Connections between non-trivial Berry's phase and topology change are discussed. \end{abstract}
\maketitle



\section{Introduction}

Differential operators on metric graphs provide an excellent opportunity to 
study connections between topology and spectral theory using differential operators
having
eigenfunctions that often can be calculated explicitly despite rather complicated
topological structure of the system. In these systems the edges determine only the geometric properties,
while topological structure is given by the
vertex conditions. Standard vertex conditions (continuity 
of the function and Kirchhoff condition on its first derivatives) used
in many examples are uniquely determined by the metric graph and give no possibility to
change the topology of the system. Using more sophisticated vertex conditions, still ensuring that
the operator is self-adjoint, allows one to consider families of systems exhibiting different topologies 
while using the same set of edges, ({\it i.e.} the same metric structure) (see Section \ref{SecS}, where we describe
these conditions in detail). 

Our attention was drawn to the subject by
\cite{Shapere2012}, where relations between vertex conditions connecting end points in two intervals
and possible topology of the corresponding graph are thoroughly discussed (following pioneering work
 \cite{Balachandran1995}). 
 In particular, the authors were interested in Berry's phase  appearing in adiabatic evolution of quantum systems (see \cite{Berry1984,SchWi1}).
 In particular the authors raised the question
(see \cite[middle of page 4]{Shapere2012}):
\begin{quotation}
whether  a non-trivial
geometric (Berry) phase can be observed, despite bases of real
eigenfunctions can be used throughout?
\end{quotation}
In the current note we answer this question positively by considering a periodic
family of metric graphs (similar to one studied in \cite{Shapere2012})
exhibiting a non-trivial geometric phase equal to $ \pi $ 
even though all eigenfunctions can be chosen real-valued: the normalised eigenfunctions of the system
are multiplied by $ -1= e^{i \pi} $ after passing one period when the system returns back to the original state.
Our model can be seen as a figure eight graph where the central vertex connects four
end points. The chosen family of vertex conditions depends on just one parameter $ \theta $ and is 
periodic with period $ 2 \pi$. For certain special values of the parameter the vertex conditions do not
properly connect all four end points and the graph turns  into either two independent cycles or just one cycle.
The suggested family of vertex conditions is also real, implying that the eigenfunctions can be chosen real-valued.
This property puts a rather strong restriction on the
system making  the existence of a non-trivial geometric phase even more striking: once
eigenfunctions are fixed for one value of the
parameter, the normalised eigenfunctions for
all other values are uniquely determined by
continuity.
Note that, without real-valued
eigenfunctions, Berry’s phase was observed
 for an  even simpler system of
just two half-lines coupled together at one degree two vertex \cite{ExGr}.

The article is organised as follows. In Section \ref{SecS} we revise the notion of metric graph
and describe relations between vertex conditions and topology of the graph. Our explicit model is presented
in Section \ref{SecModel}. Spatial symmetry of the model described in Section \ref{SecSymmetry}
allowing to treat even and odd eigenfunctions separately. The spectrum is determined in Section \ref{SecSpectrum};
it is independent of the parameter in the special case of equal edge lengths. The eigenfunctions are calculated
explicitly depending continuously on the parameter. It appears that after one period the eigenfunctions
are multiplied by $ -1 $ generating a geometric phase equal to $ \pi: e^{i \pi} = -1$.

More details can be found in Master's theses \cite{Tibbling2025,Shubin}.

\section{Topology of metric graphs and vertex conditions} \label{SecS}

In our construction we are going to follow the approach to spectral theory of metric graphs
described in detail in \cite{Book}, see also \cite{BeKu,Mugnolo}.
Let $ \Gamma $ be a metric graph formed by a finite number of {\it edges} $ E_n = [x_{2n-1}, x_{2n}], \; n = 1,2, \dots, N, $
each identified with a compact interval on an individual copy of the real line $ \mathbb R$.
Then the vertices in $ \Gamma $ can be seen as a partition of the set of end points $ \mathbf V = \{ x_j \}_{j=1}^{2N} $
into equivalence classes $ V_m, \; m =1,2, \dots, M$ ($ V_i \cap V_j = \emptyset, i \neq j, \bigcup_{m=1}^M V_m = \mathbf V$),
called {\it vertices}.
The differential operator, say the Laplacian $ - \frac{d^2}{dx^2} $, is well-defined acting on the functions
determined on the edges. The functions on different edges are independent and therefore
the differential operator is not symmetric. To make the operator self-adjoint one introduces vertex conditions
connecting limiting values of the functions 
$$ u(x_j) := \lim_{x \rightarrow x_j} u(x) $$
and their normal derivatives
$$ \partial u(x_j) = (-1)^{j+1} \lim_{x \rightarrow x_j} u'(x_j), $$
where the limits are taken from inside the intervals. 
 The derivatives are taken in the direction pointing inside the edges, hence the
extra sign in the formula. For each vertex $ V_m $ one introduces the vectors
$$ \vec{U}_m := \begin{pmatrix}
u(x_{i_1}) \\
u(x_{i_2}) \\
\vdots \\
u(x_{i_{d_m}})
\end{pmatrix}_{x_{i_j} \in V_m}, \quad \partial  \vec{U}_m := \begin{pmatrix}
\partial u(x_{i_1}) \\
\partial u(x_{i_2}) \\
\vdots \\
\partial u(x_{i_{d_m}})
\end{pmatrix}_{x_{i_j} \in V_m} . $$
The dimension of the vectors coincides with the degree $d_m $
of the vertex $ V_m $ -- the number of end points joined at $ V_m $.
Then the most general vertex conditions at $ V_m $ can be written as (see (3.21) in \cite{Book})
\begin{equation}
i (\mathbf S_m - \mathbf I) \vec{U}_m = (\mathbf S_m + \mathbf I) \partial \vec{U}_m,
\end{equation}
where $ \mathbf S_m $ is any unitary $ d \times d $ matrix.
This parametrisation of all possible vertex conditions grew up from Kostrykin-Schrader parametrisation
using pairs of matrices \cite{KoSch} and first appeared in \cite{Ha1,Ha2,KuNo}.

If the matrix $ \mathbf S_m $ has block-diagonal form, then the vertex $ V_m $ can be divided into two
(or more) smaller vertices so that the vertex conditions are introduced only joining limiting values
belonging to the new vertices. Therefore to have vertex conditions that correctly reflect the topology of the
graph the matrices $ \mathbf S_m $ should be chosen not only {\bf unitary}, but also being {\bf not block-diagonal}.

To get metric graph models with changing topology one may simply consider families of  unitary matrices
$ \mathbf S_m = \mathbf S_m (\theta) $ depending on a parameter $ \theta $ and having block-diagonal structure
for certain values of the parameter. This would correspond to splitting of the vertex for these values of $ \theta $,
while preserving original topological structure for all other values of the parameter. 

The dependence between vertex conditions and topology of the metric graph was first discussed in
\cite{KuSt} and is well-described in Section 3.3.3 of \cite{Book}. Note that vertex conditions can be used
not only to reflect different topologies but also to model physical properties of the system since 
the parameter $ \mathbf S_m $ can be seen as the vertex scattering matrix at unit energy
\cite{AsKuUs,ExTa,ExTa2,KuMa,KuEn,BaEx,KuSe}.

\section{The model} \label{SecModel}

A differential operator on metric graphs can be seen as a triple consisting of a metric graph, a differential operator
acting on the edges and vertex conditions.  The role of vertex conditions is two-fold: they are needed to make 
the differential operator self-adjoint and they determine how different edges are attached to each other, {\it i.e.} the topology of the metric graph. 
In what follows we consider
metric graphs on two edges of lengths $ \ell_1 $ and $\ell_2$
$$ E_1 = [x_1, x_2] \equiv  [-\ell_1/2,\ell_1/2] \quad  \mbox{and}  \quad E_2 = [x_3, x_4] \equiv [-\ell_2/2,\ell_2/2]. $$
The differential operator acting in the Hilbert space $ L_2 (\Gamma) = L_2 (E_1) \oplus L_2 (E_2) $ will
always be the Laplacian
\begin{equation} \label{tau}
\tau = - \frac{d^2}{dx^2}. 
\end{equation}
The vertex conditions will be given by the following one-parameter family of unitary and Hermitian matrices:
\begin{equation}
  \mathbf S_\theta = \begin{pmatrix}
    0&\sin\theta&0&\cos\theta\\
    \sin\theta&0&\cos\theta&0\\
    0&\cos\theta&0&-\sin\theta\\
    \cos\theta&0&-\sin\theta&0
  \end{pmatrix},  \quad \quad \theta \in [0,2 \pi],
\end{equation}
via the formula
\begin{equation}\label{eq_bc_ss}
  i(\mathbf S_\theta -I)\vec{u} = (\mathbf S_\theta +I)\partial {\vec{u}},
\end{equation}
where the boundary values of the functions involve their values at the end points of the intervals as well as
their normal derivatives:
\begin{equation}\label{eq:limit.values}
  \vec{u} = \begin{pmatrix}u (x_1)\\u (x_2)\\u (x_3)\\u (x_4)\end{pmatrix},\quad
  \partial {\vec{u}} = \begin{pmatrix}\partial {u}(x_1)\\\partial {u}(x_2)\\\partial {u}(x_3)\\\partial {u}(x_4)\end{pmatrix}
  :=  \begin{pmatrix}u' (x_1)\\-u' (x_2)\\u' (x_3)\\- u' (x_4)\end{pmatrix} .
\end{equation}
Since the vertex conditions are scaling-invariant, {\it i.e.} the matrix $ \mathbf S_\theta $ is not only unitary, but also Hermitian,
the ranges of the matrices $ \mathbf S_\theta - \mathbf I $  and $ \mathbf S_\theta + \mathbf I $ are orthogonal to each other (see formula (3.33) in \cite{Book}) and therefore
the vertex conditions can be explicitly written as follows, separating limiting values of the functions from the normal derivatives:
\begin{equation} \label{vc}
\begin{array}{l}
\displaystyle \begin{pmatrix}
-1 & \sin \theta & 0 & \cos \theta \\
\sin \theta & -1 &  \cos \theta & 0 \\
0 & \cos \theta & -1 & - \sin \theta \\
\cos \theta & 0 & - \sin \theta & -1
\end{pmatrix}  \begin{pmatrix}u (x_1)\\u (x_2)\\u (x_3)\\u (x_4)\end{pmatrix}  = \vec{0}, \\[10mm]
\begin{pmatrix}
1 & \sin \theta & 0 & \cos \theta \\
\sin \theta & 1 &  \cos \theta & 0 \\
0 & \cos \theta & 1 & - \sin \theta \\
\cos \theta & 0 & - \sin \theta & 1
\end{pmatrix}  \begin{pmatrix}u' (x_1)\\-u' (x_2)\\u' (x_3)\\- u' (x_4)\end{pmatrix} = \vec{0}.
\end{array}
\end{equation}
Note that the matrices above have rank $2$ for any value of $ \theta$; hence we shall always have two conditions on function values and two conditions on
the derivatives.

The corresponding operator acting on the functions from the Sobolev space $ W_2^2 (E_1) \oplus W_2^2 (E_2) $ satisfying vertex conditions \eqref{vc}
will be denoted by $ L^\theta $ and it is self-adjoint for any value of $ \theta \in [0,2 \pi].$

Note that since the matrix $ \mathbf S_\theta $ is not only unitary, but also Hermitian, the corresponding vertex scattering matrix $ S_{\mathbf v} $
does not depend on the energy parameter $k, \; k^2 = \lambda$, and is given by
\begin{equation}
\mathbf S_{\mathbf v} (k) \equiv \mathbf S_\theta.
\end{equation}

This form of vertex conditions was first used in \cite{KuEn} to model ballistic scattering of electrons 
on a wire coupled to a ring with magnetic flux. Similar vertex conditions determined by the matrix
$$ \begin{pmatrix}
0 & 0 & \sin \theta & \cos \theta \\
0 & 0 & - \cos \theta & \sin \theta \\
\sin \theta & - \cos \theta & 0 & 0 \\
\cos \theta & \sin \theta & 0 & 0 
\end{pmatrix} $$
instead of $ \mathbf S_\theta $ was used in \cite{KuSe}, where model with anomalous dependence
of the spectrum on the magnetic fluxes was considered.

\begin{figure}[H] 
\centering
\begin{tikzpicture}[scale=0.6]
  \draw (-1.5,0) circle (1.5cm);
  \draw (2,0) circle (2cm);
  \filldraw (0,0) circle (3pt);
  \node at (0.4,0.35) {$x_4$};
  \node at (-0.45,0.35) {$x_1$};
  \node at (0.4,-0.35) {$x_3$};
  \node at (-0.4,-0.35) {$x_2$};
\end{tikzpicture}
\caption{Edges parametrised.}
\label{fig:edges.param}
\end{figure}
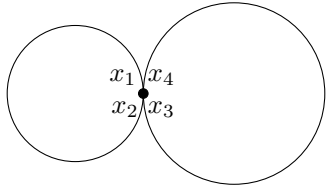

Depending on $ \theta $,  metric graphs with different topological structure correspond to the operator $ L^\theta$.
Thus, for all $ \theta \neq 0, \frac{1}{2} \pi, \pi, \frac{3}{2} \pi, 2 \pi $ all the end points $ x_1, x_2, x_3, $ and $ x_4 $ are connected 
together by the vertex conditions and the corresponding metric graph is the figure eight graph depicted in Fig. \ref{fig:edges.param} -- the matrix $ \mathbf S_\theta $ is not block-diagonal.

For other values of $ \theta $ the matrix $ \mathbf S_\theta $ is block-diagonal and
the end points can be divided into pairs so that the corresponding metric graph is either formed by two loops
of lengths $ \ell_1 $ and $ \ell_2 $ 
$$
\begin{array}{lcccc}
\displaystyle \theta = \frac{1}{2} \pi & \Rightarrow & \displaystyle  
\left\{
\begin{array}{l}
\sin \theta = 1 \\ \cos \theta = 0
\end{array} 
\right.&  \Rightarrow  &\begin{array}{l}
\left\{ \begin{array}{ccc}
u(x_1) & = & u(x_2) \\
u'(x_1) & = & u'(x_2) \end{array} \right.       \\[3mm]   \left\{ \begin{array}{ccc}
u(x_3) & = & - u(x_4) \\
u'(x_3) & = & - u'(x_4) \end{array} \right.    
\end{array}; \\[10mm]
\displaystyle \theta= \frac{3}{2} \pi & \Rightarrow & 
\left\{ \begin{array}{l}
\sin \theta = -1 \\
 \cos \theta = 0 \end{array}  \right.
 &  \Rightarrow  & \begin{array}{l}
 \left\{ \begin{array}{ccc}
u(x_1) & = & - u(x_2) \\
u'(x_1) & = & - u'(x_2) \end{array} \right.       \\[3mm]
  \left\{ \begin{array}{ccc}
u(x_3) & = & u(x_4) \\
u'(x_3) & = &  u'(x_4) \end{array} \right.   
\end{array};
\end{array} $$
or by single loop of length $ \ell_1 + \ell_2 $ with reflectionless conditions at two internal points:
$$
\begin{array}{lcccc}
\displaystyle \theta = 0 & \Rightarrow & \displaystyle  
\left\{
\begin{array}{l}
\sin \theta = 0 \\ \cos \theta = 1
\end{array} 
\right.&  \Rightarrow  &\begin{array}{l}
\left\{ \begin{array}{ccc}
u(x_1) & = & u(x_4) \\
u'(x_1) & = & u'(x_4) \end{array} \right.       \\[3mm]   \left\{ \begin{array}{ccc}
u(x_2) & = &  u(x_3) \\
u'(x_2) & = & u'(x_3) \end{array} \right.    
\end{array}; \\[10mm]
\displaystyle \theta=  \pi & \Rightarrow & 
\left\{ \begin{array}{l}
\sin \theta = 0 \\
 \cos \theta = -1 \end{array}  \right.
 &  \Rightarrow  & \begin{array}{l}
 \left\{ \begin{array}{ccc}
u(x_1) & = & - u(x_4) \\
u'(x_1) & = & - u'(x_4) \end{array} \right.       \\[3mm]
  \left\{ \begin{array}{ccc}
u(x_2) & = & - u(x_3) \\
u'(x_2) & = &  - u'(x_3) \end{array} \right.   
\end{array}.
\end{array} $$

\begin{figure}[H]
\centering
\begin{tikzpicture}[scale=0.8]
  \draw (2.5,1.5) circle (0.5cm);
  \draw (3.5,1.5) circle (0.5cm);
  \filldraw (3,1.5) circle (0.10cm);
  \draw (0.6,2) circle (0.5cm);
  \draw (-0.6,2) circle (0.5cm);
  \draw (0.1,2) circle (0.07cm);
  \draw (-2.5,1.5) circle (0.5cm);
  \draw (-3.5,1.5) circle (0.5cm);
   \filldraw (-3,1.5) circle (0.10cm);
  \draw (-5,0) circle (0.7cm);
  \draw (-5,-0.7) circle (0.07cm);
  \draw (-5,0.7) circle (0.07cm);
  \draw (-2.5,-1.5) circle (0.5cm);
  \draw (-3.5,-1.5) circle (0.5cm);
   \filldraw (-3,-1.5) circle (0.10cm);
  \draw (0.6,-2) circle (0.5cm);
  \draw (-0.6,-2) circle (0.5cm);
  \draw (-0.1,-2) circle (0.07cm);
  \draw (2.5,-1.5) circle (0.5cm);
  \draw (3.5,-1.5) circle (0.5cm);
   \filldraw (3,-1.5) circle (0.10cm);
  \draw (5,0) circle (0.7cm);
  \draw[-{Straight Barb}] (3,2.5) .. controls (0,3.2) .. (-3,2.5);
  \draw[-{Straight Barb}] (-3,-2.5) .. controls (0,-3.2) .. (3,-2.5);
  \draw[-{Straight Barb}] (-5.5,1.5) .. controls (-6,0) .. (-5.5,-1.5);
  \draw[-{Straight Barb}] (5.5,-1.5) .. controls (6,0) .. (5.5,1.5);
  \node at (6.7,0) {$\theta=0$};
  \node at (-6.7,0) {$\theta=\pi$};
  \node at (0,3.7) {$\theta = \pi/2$};
  \node at (0,-3.7) {$\theta = 3\pi/2$};
\end{tikzpicture}
\caption{Cycle of topological changes corresponding to $\mathbf S_\theta$.}
\label{fig:topology_change}
\end{figure}
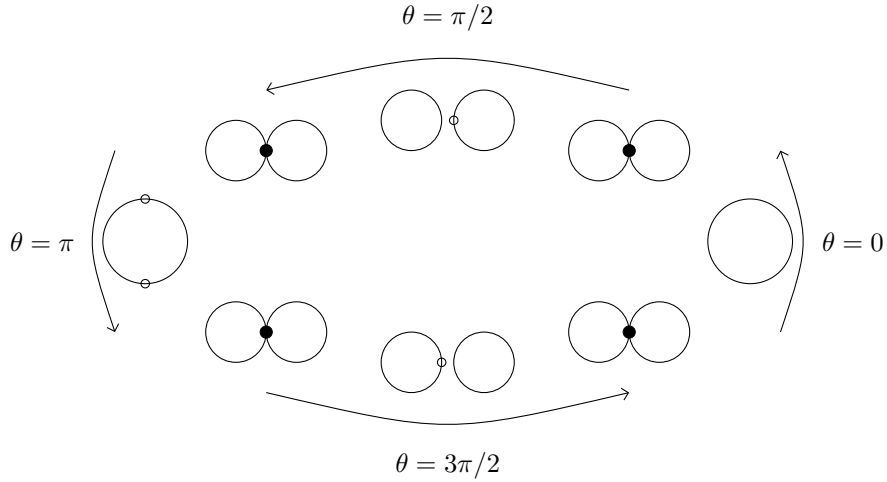

This dependence of the topology of the metric graph $ \Gamma $ on the parameter $ \theta $ is best illustrated by Fig. \ref{fig:topology_change}.
In this figure we use the following notations for the vertices:
\begin{itemize}
\item filled large circles: degree four vertex with  non-separable conditions involving parameter $ \theta \neq 0, \frac{1}{2} \pi, \pi, \frac{3}{2} \pi, 2 \pi$;
\item empty small circles: degree two vertices with reflectionless conditions involving multiplication of the function and its first derivative by $-1$;
\item no circle:  degree two vertices with standard conditions (continuity of the function and its first derivative).
\end{itemize}

We see that the topology of the metric graph associated with the operator $ L^\theta $ changes with $ \theta$ and 
the corresponding graphs are not even connected for $ \theta =  \frac{1}{2} \pi, \frac{3}{2} \pi $.
The number of independent cycles also changes and is equal to $1$ for $ \theta = 0, \pi$.

\section{Graph's symmetry} \label{SecSymmetry}

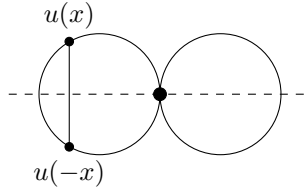
\begin{figure}[H]
\centering
\begin{tikzpicture}[scale=0.8]
  \draw (-1,0) circle (1cm);
  \draw (1,0) circle (1cm);
  \filldraw (0,0) circle (3pt);
  \filldraw (-1.5,-0.87) circle (2pt);
  \filldraw (-1.5,0.87) circle (2pt);
  \draw[dashed] (-2.5,0) -- (2.5,0);
  \node at (-1.5,-1.3) {$u(-x)$};
  \node at (-1.5,1.3) {$u(x)$};
  \draw (-1.5,-0.87) -- (-1.5,0.87);
\end{tikzpicture}
\caption{Horizontal symmetry.} \label{fig:sym}
\end{figure}

The constructed operator possesses the following symmetry, which will help us to calculate its eigenfunctions.
Let us denote by $ J $ the unitary operator corresponding to the horizontal symmetry of the figure eight graph (see Fig. \ref{fig:sym}):
\begin{equation}
\begin{array}{cccc}
\displaystyle J: & \displaystyle L_2(E_1) \oplus L_2 (E_2) & \rightarrow & \displaystyle  L_2(E_1) \oplus L_2 (E_2);  \\[3mm]
&  \displaystyle  (u_{\rm 1} (x), u_{\rm 2} (x) ) & \mapsto &   \displaystyle  (u_{\rm 1} (-x), u_{\rm 2} (-x) ). 
\end{array}
\end{equation}

The symmetry transformation  $J $ acts on the vectors  $ \vec{u} $ and $ \partial \vec{u}$ of limiting values introduced in \eqref{eq:limit.values} 
as multiplication by the matrix
$$ \mathbf J :=   \begin{pmatrix}
0 & 1 & 0 & 0 \\
1 & 0 & 0 & 0 \\
0 & 0 & 0 & 1 \\
0 & 0 & 1 & 0
\end{pmatrix}.
$$
We have obviously
$$ \mathbf J \mathbf S_\theta = \mathbf S_\theta \mathbf J,$$
which together with
$$ J \tau = \tau J ,$$
(where $ \tau $ was defined in \eqref{tau})
implies that the operator $ L^\theta $ commutes with the symmetry operator $ J $
\begin{equation}
J L^\theta = L^\theta J. 
\end{equation}  

It follows that the eigenfunctions may be divided into two classes to be studied separately:
\begin{itemize}
\item even eigenfunctions satisfying $ J \psi= \psi $;
\item odd eigenfunctions satisfying  $ J \psi = - \psi$.
\end{itemize}

The defined operator is real in the sense that if $ u $ belongs to the domain of the operator, then the complex conjugate $ \overline{u} $ also
belongs to the domain and $ \tau \overline{u} = \overline{\tau u}. $
Therefore if $ \psi_\lambda $ is an eigenfunction corresponding to the eigenvalue $ \lambda$, then $ \overline{\psi}_\lambda $ is
also an eigenfunction corresponding to the same eigenvalue (which is real since the operator is self-adjoint). Hence both $ \Re \psi := \frac{\psi + \overline{\psi}}{2} $
and $ \Im \psi := \frac{\psi - \overline{\psi}}{2i} $ are eigenfunctions implying that without loss of generality the eigenfunctions may
always be chosen not only even/odd, but also real-valued.

\section{The secular equation}\label{sec:secular.equation} \label{SecSpectrum}

\subsection*{Non-zero eigenvalues}

To get the secular equation describing the non-zero spectrum of the operator we shall use secular polynomials
well described in Section 6.1 of \cite{Book}. We introduce  the secular polynomial $ P(z_1, z_2) $ determined by
\begin{equation}
\begin{array}{ccl}
\displaystyle P(z_1, z_2) & = &  \displaystyle   \det \Big(\mathbf S_{\rm e} - \underbrace{\mathbf S_{\rm v}}_{\displaystyle  
 =  \mathbf S_\theta} \Big) \\
& = & \displaystyle  \det \left(\begin{pmatrix}
0 & z_1 & 0 & 0 \\
z_1 & 0 & 0 & \\
0 & 0 & 0 & z_2 \\
0 & 0 & z_2 &0
\end{pmatrix} 
-
\begin{pmatrix}
    0&\sin\theta&0&\cos\theta\\
    \sin\theta&0&\cos\theta&0\\
    0&\cos\theta&0&-\sin\theta\\
    \cos\theta&0&-\sin\theta&0
  \end{pmatrix}
  \right) \\[10mm]
  & = &  \displaystyle 
  \left( (z_1 - \sin \theta) (z_2 + \sin \theta) - \cos^2 \theta \right)^2.
  \end{array}
\end{equation}
Then the spectrum of the operator is given by zeroes of the trigonometric polynomial
\begin{equation}
\begin{array}{ccl}
p(k) & := & P(e^{ik \ell_1}, e^{i k \ell_2}) \\
& = & \displaystyle - 4 e^{i k (\ell_1+\ell_2)/2}  \left( \sin \left(k \frac{\ell_1 + \ell_2}{2}\right) + \sin \theta \cdot  \sin \left( k \frac{\ell_1-\ell_2}{2} \right) \right)^2.
\end{array}
\end{equation}
Ignoring the non-vanishing factor we conclude that all non-zero eigenvalues are double degenerate and are given by the solutions  of
the equation
\begin{equation}
\label{secular}
 \sin k \frac{\ell_1 + \ell_2}{2} + \sin \theta \cdot  \sin k \frac{\ell_1-\ell_2}{2} = 0. 
 \end{equation}
 It will be proven in the next section that for each zero one eigenfunction can be chosen even and one odd.
 The secular equation \eqref{secular} can also be derived directly by noting that any even eigenfunction on $ \Gamma $ is given by
 $$  \psi_n  (x) = a_j \cos n \pi x, \quad x \in E_j, j= 1,2 $$
 and every odd by
 $$ \psi_n  (x) = b_j \sin n \pi x, \quad x \in E_j, j= 1,2.$$

In the special case $ \ell_1 = \ell_2 = \ell $ the spectrum is independent of $ \theta $ and is given by 
\begin{equation}
\sin k \ell = 0 \quad \Rightarrow \quad  k_n = \frac{ \pi}{\ell} n , \quad n = 1,2,3, \dots
\end{equation}

\subsection*{The eigenvalue $ \lambda = 0 $} To determine the multiplicity of the zero eigenvalue one needs to 
repeat our
analysis
 taking into account that the solution of the eigenfunction equation  $ \psi \psi (x) = 0 $
 on the edges is given by a linear function instead
of the exponentials. It is easier to consider possible even and odd eigenfunctions separately.

\subsubsection*{Even eigenfunctions}
Every even eigenfunction $ \psi $ corresponding to $ \lambda = 0 $ is equal to a constant function on each of
the edges:
$$ \psi (x) =  a_j \quad x \in E_j, j= 1,2. 
$$
Its boundary values are
$$ \vec{\psi} = \begin{pmatrix}
a_1 \\
a_1 \\
a_2 \\
a_2
\end{pmatrix}, \quad \partial \vec{\psi} = \vec{0}.  $$
Substitution into the vertex conditions \eqref{vc} leads to just two linear equations for $ a_1 $ and $a_2 $
$$
\left\{
\begin{array}{rcrcl}
(-1 + \sin \theta) \; a_1 & + & \cos  \theta \; a_2 & = & 0, \\
\cos \theta \; a_1 & - & (1+ \sin \theta) \; a_2 & = & 0.
\end{array} \right.
$$
The determinant of the linear system is always equal to zero: $ 1 - \sin^2 \theta - \cos^2 \theta = 0 $,
and all four entries do not vanish simultaneously. Hence there is precisely one even eigenfunction corresponding to $ \lambda = 0 $  for any value of $ \theta,$
independently of whether the corresponding metric graph is connected or not.

\subsubsection*{Odd eigenfunctions}
Every odd eigenfunction corresponding to $ \lambda = 0 $ is equal to a linear function on each of
the edges:
$$ \psi (x) = b_j x, \quad x \in E_j, j= 1,2.
$$
We calculate its limit values at the end points:
$$ \vec{\psi} = \begin{pmatrix}
- b_1 \ell_1/2  \\
b_1 \ell_1/2  \\
- b_2 \ell_2/2  \\
b_2 \ell_2/2 
\end{pmatrix}, \quad \partial \vec{\psi} =  \begin{pmatrix}
b_1   \\
- b_1   \\
 b_2   \\
-b_2  
\end{pmatrix}.  $$
Then vertex conditions  \eqref{vc} imply the following linear system after excluding identical equations
$$
\begin{pmatrix}
\frac{\ell_1}{2} (1+ \sin \theta)  &  \frac{\ell_2}{2} \cos \theta  \\
\frac{\ell_1}{2} \cos \theta & \frac{\ell_2}{2}  (1-\sin \theta)\\
1- \sin \theta & - \cos \theta \\
- \cos \theta & 1+ \sin \theta
\end{pmatrix}
\begin{pmatrix} b_1 \\
b_2
\end{pmatrix} = \vec{0}.
$$
The matrix has rank $2$ for any value of $ \theta$, hence the system has only  the trivial solution
implying that no odd eigenfunctions are present.

We conclude that $ \lambda _0= 0 $ is a simple eigenvalue and the corresponding eigenfunction is
even.

By combining all the results of this section, we obtain the following theorem.

\begin{theorem}\label{secular.theorem}
The spectrum of the Laplacian $ L^\theta $ defined on the domain of functions from 
the Sobolev space $ W_2^2 (E_1) \oplus W_2^2 (E_2) $ satisfying vertex conditions \eqref{vc}
is described below:
\begin{itemize}
\item the ground state $ \lambda_0 = 0 $ is a simple eigenvalue and the corresponding eigenfunction is even;
\item all non-zero eigenvalues $ \lambda_n = k_n^2 $ are double degenerate and are given by the zeroes of
the trigonometric polynomial
\begin{equation}
 \sin \left( k \frac{\ell_1 + \ell_2}{2} \right) + \sin \theta \cdot  \sin \left( k \frac{\ell_1-\ell_2}{2}  \right).
 \end{equation}
\end{itemize}
In the special case $ \ell_1 = \ell_2 = :\ell$ the spectrum is independent of $ \theta $ and is given by:
\begin{equation}
\sigma (L^\theta) \vert_{\ell_1 = \ell_2} = \left\{ 0,  \left(\frac{\pi}{\ell} \right)^2,  \left(\frac{\pi}{\ell} \right)^2,  4 \left(\frac{\pi}{\ell} \right)^2 ,  4 \left(\frac{\pi}{\ell} \right)^2 , \dots   \right\}
\end{equation}
\end{theorem}

The spectrum can be explicitly calculated not only for $ \ell_1 = \ell_2 $, but for the special values of $ \theta $:
\begin{itemize}
\item $ \theta = 0, \pi $, then the spectrum is:
$$  \begin{array}{ll}
0, & \mbox{simple eigenvalue}, \\
n^2 \left(\frac{2\pi}{\ell_1+\ell_2} \right)^2 , \;  n = 1,2,3,  \dots, &
\mbox{double eigenvalues.}
\end{array} $$ 
\item $ \theta = \pi/2 $, then the spectrum is:
$$  \begin{array}{ll}
0, & \mbox{simple eigenvalue}, \\
(2n)^2 \left(\frac{\pi}{\ell_1} \right)^2 , \;  n = 1,2, 3,\dots, &
\mbox{double eigenvalues,} \\
(2n+1)^2 \left(\frac{\pi}{\ell_2} \right)^2 , \; n = 0,1,2, \dots, &
\mbox{double eigenvalues.}
\end{array} $$ 
\item $ \theta = 3 \pi/2 $, then the spectrum is: 
$$  \begin{array}{ll}
0, & \mbox{simple eigenvalue}, \\
(2n+1)^2 \left(\frac{\pi}{\ell_1} \right)^2 , \; n = 0, 1,2, \dots, &
\mbox{double eigenvalues,} \\
(2n)^2 \left(\frac{\pi}{\ell_2} \right)^2 , \; n = 1, 2, 3, \dots, &
\mbox{double eigenvalues.}
\end{array} $$

\end{itemize}

\section{The eigenfunctions}

In this section we are going to calculate the eigenfunctions  denoted by $ \psi $.
We have to bear in mind that 
\begin{itemize}
\item the operator is real, {\it i.e.} invariant under 
complex conjugation, and therefore the eigenfunctions can be chosen real-valued,
\item the operator is symmetric $  J L^\theta = L^\theta J $ and therefore
one may look separately for even and odd eigenfunctions.
\end{itemize}
To make our presentation even more transparent we consider the case of equal unit lengths
$ \ell_1 = \ell_2 = 1 $, in which case all operators $ L^\theta $ are isospectral to each other.

\subsection*{Even eigenfunctions} Every even eigenfunction is given by
$$
\psi_n  (x) = a_j \cos n \pi x, \quad x \in E_j, j= 1,2.
$$
We just need to determine the real constants $ a_1 $ and $ a_2 $ so that $ \psi_n $ 
is normalised and depends continuously on $ \theta$.

The first condition comes from the normalisation of the eigenfunction
$$1 = \parallel  \psi_n \parallel^2 = (a_1^2 + a_2^2 ) \int_{-1/2}^{1/2} \cos^2 n \pi x \; dx = \frac{1}{2} (a_1^2 + a_2^2) \Rightarrow $$
\begin{equation} \label{eq2}
a_1^2 + a_2^2 = 2.
\end{equation}

Checking that the function satisfies the vertex condition we shall need to separate the cases when $ n $ is even and odd.

\subsubsection*{Case $ n = 2 m +1 , \quad m = 0,1,2, \dots $} $ $ \newline

\noindent
Substituting the limiting values of the function and its first derivative into the vertex conditions \eqref{vc}
we get:
$$
\vec{\psi}_n = \vec{0}, \quad   \partial  \vec{\psi}_n =  (-1)^m (2m +1) \pi \begin{pmatrix}
a_1 \\
a_1 \\
a_2 \\
a_2
\end{pmatrix}
\Rightarrow
$$
\begin{equation} \label{eq1}
(1+ \sin \theta) \; a_1 + \cos \theta \; a_2 = 0.
\end{equation}

Our task now is to determine continuous  functions $ a_1 (\theta) $ and $ a_2 (\theta) $ satisfying equations
\eqref{eq2} and \eqref{eq1}. Let us fix a normalisation  compatible with \eqref{eq2} and \eqref{eq1} by assuming that
\begin{equation} \label{eq3}
a_1(0) = 1, \quad a_2 (0) =  - 1.
\end{equation}
Excluding $ a_2 $ from equation \eqref{eq1} using equation \eqref{eq2} we get
\begin{equation} \label{eq11}
 a_1 (\theta) = \pm \frac{\cos \theta}{\sqrt{1 + \sin \theta}}. 
 \end{equation}
The corresponding $ a_2 (\theta ) $ is then given by
\begin{equation} \label{eq12}
 a_2 (\theta) = \mp \sqrt{1  + \sin \theta}. 
 \end{equation}
To ensure continuous
dependence of $ a_j (\theta) $ on $ \theta $ we need to introduce the sign function
\begin{equation}
\sigma (\theta) = \left\{
\begin{array}{ll}
1, & 0 \leq \theta < 3 \pi/2, \\
-1 &  3 \pi/2 \leq \theta \leq 2 \pi.
\end{array} \right. 
\end{equation}
Then the amplitudes given by the formulas below depend continuously on $ \theta $, the function satisfies the vertex conditions and is normalised
\begin{equation}
a_1 (\theta) = \sigma (\theta) \frac{\cos \theta}{\sqrt{1 +\sin \theta}}, \quad a_2 (\theta) = - \sigma(\theta) \sqrt{1 + \sin \theta}.
\end{equation}
These formulas can be simplified using trigonometric identities
$$
\begin{array}{l}
\displaystyle  \sqrt{1 + \sin \theta} = \sqrt{\sin^2 \theta/2 + \cos^2 \theta/2 + 2 \sin \theta/2 \cos \theta/2} = \vert \sin \theta/2 + \cos \theta/2 \vert; \\[3mm]
\displaystyle \sigma_1 (\theta) = {\rm sgn}\; \Big( \sin \theta/2 + \cos \theta/2 \Big), \quad \quad \theta \in [0,2 \pi];
\end{array}$$
as follows
\begin{equation} \label{eq5}
a_1 (\theta) = \cos \frac{\theta}{2} - \sin \frac{\theta}{2}, \quad \quad a_2 (\theta) =   - \cos \frac{\theta}{2} - \sin \frac{\theta}{2}.
\end{equation}
It would have been possible to derive formulas \eqref{eq5} directly from \eqref{eq11} and \eqref{eq12}\footnote{We follow this path calculating below the amplitudes for even values of $ n$.} and choosing
$$ a_2 (\theta) = - \cos \frac{\theta}{2} - \sin \frac{\theta}{2}, $$
which in turn implies
$$ 
\begin{array}{ccl}
a_1 (\theta) & = & \displaystyle - \frac{\cos \theta}{1+ \sin \theta} a_2 (\theta) = \frac{\cos \theta}{1 + \sin \theta} \Big( \cos \theta/2 + \sin \theta/2 \Big)  \\[3mm]
& = &  \displaystyle \frac{\cos^2 \theta/2 - \sin^2 \theta/2}{(\cos \theta/2 + \sin \theta/2)^2}  \Big( \cos \theta/2 + \sin \theta/2 \Big)  \\[5mm]
& = & \displaystyle \cos \theta/2 - \sin \theta/2.
\end{array}
$$
The role of the sign function is hidden in this approach.

\subsubsection*{Case $ n = 2 m , \quad m = 1,2, 3,  \dots $} $ $ \newline

\noindent The limiting values are given by
$$
\vec{\psi}_n = (-1)^m \begin{pmatrix}
a_1 \\
a_1 \\
a_2 \\
a_2
\end{pmatrix} \quad \partial  \vec{\psi}_n  = \vec{0}, 
$$
and lead to the equation
\begin{equation} \label{eq13}
(-1+ \sin \theta) \; a_1 + \cos \theta \; a_2 = 0,
\end{equation}
which should be solved together with the normalisation condition \eqref{eq2}.

This time the normalisation, compatible with \eqref{eq2} and \eqref{eq13},  can be
chosen to be given by
\begin{equation} \label{eq3}
a_1(0) = 1, \quad a_2 (0) = 1.
\end{equation}
Excluding $ a_1 $ using \eqref{eq13} and substituting into \eqref{eq2} we get
$$ a_2 (\theta) = \pm \sqrt{1 - \sin \theta} = \pm \Big( \cos \theta/2 - \sin \theta/2 \Big). $$
Then normalisation \eqref{eq3} requires
\begin{equation} \label{eq41}
 a_2 (\theta) = \cos \theta/2 - \sin \theta/2, \qquad 0 \leq  \theta  \leq \pi/2. 
 \end{equation}
We calculate also the first amplitude
\begin{equation} \label{eq42}
\begin{array}{ccl}
a_1 (\theta) & = & \displaystyle 
 \frac{\cos \theta}{1 - \sin \theta} \Big( \cos \theta/2 - \sin \theta/2 \Big) = \frac{\cos^2 \theta/2 - \sin^2 \theta/2}{\Big( \cos \theta/2 - \sin \theta/2 \Big)^2} \Big( \cos \theta/2 - \sin \theta/2
\Big) \\[5mm]
& = & \displaystyle  \cos \theta/2 + \sin \theta/2.
\end{array}
\end{equation}
Formulas \eqref{eq41} and \eqref{eq42} determine  the amplitudes  for $ 0 \leq \theta \leq \pi/2$ and there is just a unique way to extend
them for $ \pi/2 \leq \theta \leq 2 \pi$ keeping the functions continuous:
\begin{equation} \label{eq6}
a_1 (\theta) = \cos \frac{\theta}{2}  + \sin \frac{\theta}{2}, \qquad a_2 (\theta)  = \cos \frac{\theta}{2} -  \sin \frac{\theta}{2}.
\end{equation}

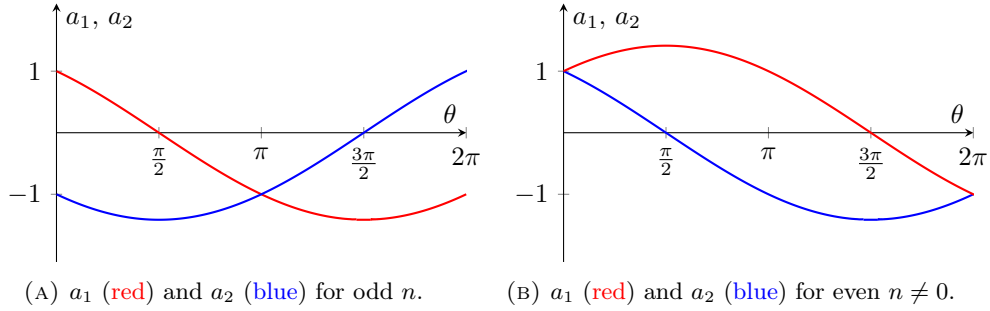
\begin{figure}[H]
\centering
\begin{subfigure}[b]{0.47\textwidth}
\centering
\begin{tikzpicture}
\begin{axis}[
  width=7cm, height=5cm,
  xlabel={$\theta$}, ylabel={$a_1,\,a_2$},
  xmin=0, xmax=6.28318,
  ymin=-2.1, ymax=2.1,
  xtick={0,1.5708,3.14159,4.71239,6.28318},
  xticklabels={$0$,$\frac{\pi}{2}$,$\pi$,$\frac{3\pi}{2}$,$2\pi$},
  ytick={-1,0,1},
  axis lines=center,
  legend pos=north west,
  legend style={font=\small},
]
\addplot[red,  thick, domain=0:4.60,       samples=300] {cos(deg(x))/sqrt(1+sin(deg(x)))};
\addplot[red,  thick, domain=4.60:4.82, samples=300] {-1.414};
\addplot[red,  thick, domain=4.82:6.28318, samples=300] {-cos(deg(x))/sqrt(1+sin(deg(x)))};
\addplot[blue, thick, domain=0:6.88318, samples=300] {-sin(deg(x)/2)-cos(deg(x)/2)};
\end{axis}
\end{tikzpicture}
\caption{$a_1 $ (\textcolor{red}{red}) and $a_2$ (\textcolor{blue}{blue}) for odd $n$.}
\end{subfigure}
\hfill
\begin{subfigure}[b]{0.47\textwidth}
\centering
\begin{tikzpicture}
\begin{axis}[
  width=7cm, height=5cm,
  xlabel={$\theta$}, ylabel={$a_1,\,a_2$},
  xmin=0, xmax=6.28318,
  ymin=-2.1, ymax=2.1,
  xtick={0,1.5708,3.14159,4.71239,6.28318},
  xticklabels={$0$,$\frac{\pi}{2}$,$\pi$,$\frac{3\pi}{2}$,$2\pi$},
  ytick={-1,0,1},
  axis lines=center,
]
\addplot[blue,  thick, domain=0:4.60,       samples=300] {cos(deg(x))/sqrt(1+sin(deg(x)))};
\addplot[blue,  thick, domain=4.60:4.82, samples=300] {-1.414};
\addplot[blue,  thick, domain=4.82:6.28318, samples=300] {-cos(deg(x))/sqrt(1+sin(deg(x)))};
\addplot[red, thick, domain=0:6.88318, samples=300] {sin(deg(x)/2)+cos(deg(x)/2)};
\end{axis}
\end{tikzpicture}
\caption{$a_1$ (\textcolor{red}{red}) and $a_2$ (\textcolor{blue}{blue}) for even $n \neq 0$.}
\end{subfigure}

\caption{Continuous eigenfunction coefficients over one cycle.}  \label{Figplot}
\end{figure}

The amplitudes depending continuously on $ \theta $ are plotted in Fig. \ref{Figplot}.
Formulas \eqref{eq5} and \eqref{eq6} determine the amplitudes as $ 4 \pi$-periodic functions of $ \theta$,
while period $ 2 \pi $ may be expected since the family $ L^\theta$ has period $ 2 \pi$.

One clearly observes that
\begin{equation} \label{Berry}
\psi_n \vert_{\theta = 2 \pi}  = \underbrace{-}_{= e^{i \pi}} \psi_n \vert_{\theta = 0}.
\end{equation}

\subsection*{Odd eigenfunctions}
The analysis is entirely analogous to the even case with the only difference that the amplitudes for even and odd values of $ n $
are exchanged. The eigenfunctions demonstrate a topological phase $ \pi $ after one period $ \theta: 0 \rightarrow 2 \pi$.

\subsection*{The ground state $ \lambda = 0 $}
The analysis is very similar to the case of even eigenfunctions $ \psi_{2m} $ with the only difference that the normalisation condition leads to
\begin{equation} \label{eq22}
a_1^2 + a_2^2 = 1,
\end{equation}
instead of \eqref{eq2}. Hence the corresponidng amplitudes are given by
 \begin{equation} \label{eq66}
a_1 (\theta) = \frac{1}{\sqrt{2}}  \Big(\cos \frac{\theta}{2}  + \sin \frac{\theta}{2} \Big) , \qquad a_2 (\theta)  =   \frac{1}{\sqrt{2}} \Big(\cos \frac{\theta}{2}  - \sin \frac{\theta}{2} \Big).
\end{equation}
The eigenfunction generates the  topological phase $ \pi $ after one period. 

\begin{theorem}
Let $ L^\theta $ be the $2 \pi$-periodic family of Laplacians on graphs on two edges determined by the vertex conditions \eqref{vc}. 
Then the eigenfunctions of the operator chosen to be  real-valued and continuous generate the geometric Berry's phase $ \pi $, {\it i.e.}
the eigenfunctions satisfy the equation
\begin{equation}
\psi_n (x) \vert_{\theta = 2 \pi} = e^{i \pi} \psi_n (x) \vert_{\theta =  0}, \quad n = 0,1,2, \dots
\end{equation}
\end{theorem}

The theorem has been proven in the special case of equal edge lengths, but the spectrum and the eigenfunctions
depend continuously on the edge lengths. Hence the geometric phase depends continuously on
the edge lengths, but it attains just two possible values of $ 0 $ and $ \pi$, hence the conclusion of the theorem holds even for not 
necessarily equal edge lengths.

\section{Conclusions and discussions}

Thus we have proven that the topological Berry phase for the considered family of graphs is non-trivial and equals $ \pi. $
The main reason for non-triviality of the geometric phase is connected with the  fact that topology of the
system changes with $ \theta $: for almost all values of $ \theta $ the graph should be seen as the figure eight graph with 
vertex conditions depending on the parameter, but for certain special values of the parameter the graph reduces either
to one or two cycles. In other words, 
not only does the number of cycles
change with $ \theta$, but also the connectivity of the
graph.
The eigenfunctions are essentially described by two amplitudes depending continuously on $ \theta $, and to gain 
a non-trivial topological phase these amplitudes have to
change sign In other words each of the amplitudes has to vanish. The corresponding eigenfunction vanishes identically
on one of the edges. In our model this appears precisely when the graph loses connectivity ($ \theta = \pi/2, 3 \pi/2$)
and reduces to two independent circles. 

Identical vanishing of eigenfunctions on some edges is a typical feature of operators on metric graphs and may happen even without 
losing connectivity. Therefore, it might be interesting to find an example of 
a family of metric
graphs with a non-trivial geometric phase whose
topology does not change.

\section*{Data availability statement}

No new data were created or analysed in this study.

\section*{Funding} 

The research of P.K. was been supported by The Swedish Research Council, grant number 2024-04650.

\section*{Author contributions}

The authors contributed to deriving, analysing the formulas and writing the manuscript. The authors read
and approved the final manuscript.

\section*{Conflict of interest}

The authors declare that there are no competing interests.

\end{document}